\documentclass[10pt]{amsart}          

\usepackage{amsmath, amscd, amsfonts, amssymb}
\usepackage{mathrsfs, amsthm, bbm}
\usepackage[all]{xy}
\usepackage{color}
\usepackage{wrapfig}		

\newtheorem{proposition}{Proposition}
\newtheorem{theorem}[proposition]{Theorem}

\theoremstyle{remark}

\newtheorem{remark}[proposition]{Remark}

\theoremstyle{definition}

\newcommand{\cl}{{\rm C\ell}}
\newcommand{\id}{{\rm id}}
\newcommand{\z}{\mathbbm Z}
\renewcommand{\r}{\mathbbm R}
\renewcommand{\c}{\mathbbm C}

\begin{document}

\title{The spinor bundle of Riemannian products}
\author{Frank Klinker}
\thanks{This text is part of the authors PhD.\ thesis \cite{Klinker2}}
\date{\today}
\address{Department of Mathematics,  University of Dortmund,  D--44221 Dortmund}
\address{{\tt frank.klinker@math.uni-dortmund.de}}
\begin{abstract} 
In this note we compare the spinor bundle of a Riemannian manifold $(M=M_1\times\cdots\times M_N,g)$ with the spinor bundles of the Riemannian factors $(M_i,g_i)$. We show that - without any holonomy conditions - the spinor bundle of $(M,g)$ for a special class of metrics is isomorphic to a bundle obtained by tensoring the spinor bundles of $(M_i,g_i)$ in an appropriate way. For $N=2$ and an one dimensional factor this construction was developed in \cite{Baum5}. Although the fact for general factors is frequently used in (at least physics) literature, a proof was missing.

I would like to thank Shahram Biglari, Mario Listing, Marc Nardmann and Hans-Bert Rademacher for helpful comments. Special thanks go to Helga Baum, who pointed out some difficulties arising in the pseudo-Riemannian case.
\end{abstract}
\maketitle

We consider a Riemannian manifold $(M=M_1\times\cdots\times M_N,g)$, which is a product of Riemannian spin manifolds $(M_i,g_i)$  and denote the projections on the respective factors by $p_i$. Furthermore the dimension of $M_i$ is $D_i$ such that the dimension of $M$ is given by $D=\sum_{i=1}^N D_i$. 

The tangent bundle of $M$ is decomposed as
\begin{equation}
T_{(x_0,\ldots,x_N)}M= p_1^*T_{x_1}M_1\oplus\cdots \oplus p_N^*T_{x_N}M_N.
\label{tangentsplitting}
\end{equation}
We omit the projections and write $TM=\bigoplus_{i=1}^{N} TM_i$. 

The metric $g$ of $M$ need not be the product metric of the metrics $g_i$ on $M_i$, but is assumed to be of the form 
\begin{multline}
g_{ab} (x)\big|_{TM_i}= {A_i}^c_a(x){g_i}_{cd}(x_i){A_i}^d_b(x),\\ 
\text{for }  D_1+\cdots +D_{i-1}+1\leq  a,b\leq D_1+\cdots+D_i, 1\leq i\leq N
\label{metricsplitting}
\end{multline}
In particular, for those metrics the splitting (\ref{tangentsplitting}) is  orthogonal, i.e.\ the frame bundle of $M$ can be reduced to a $SO(D_1)\times\cdots\times SO(D_N)$-principal bundle, and this is isomorphic to the product of the frame bundles over $M_i$. The explicit form of the isomorphism is
\begin{multline}
P_{SO}(M_1)\times P_{SO}(M_N)\ni
(E_1(x_1),\ldots ,E_N(x_N)) \\
\mapsto(A_1(x)E_1(x_1),\ldots,A_N(x)E_N(x_N))\in P_{SO}(M)
\label{isomorphism}
\end{multline}

It is clear that such a manifold need not  have a splitting of the holonomy group into subgroups of $SO(D_i)$, which would  lead to an, at least local, Riemannian product (cf. \cite[sect.\ 3.2]{JoyceBook}). Examples for such spaces are the Eguchi Hanson space -- where we have $M\subset \r\times S^3$ with metric $\alpha(r)dr^2+h_r(S^3)$ -- or warped products of metrics.

We consider spinor bundles $S_i$ over $M_i$ and we are going to construct a bundle $S$ over $M$ from these spinor bundles.
We will discuss what conditions are necessary for the bundle $S$ to be  the spinor bundle over $M$ and how we have to modify the given Clifford multiplication  on $S_i$  to extend it to $S$.  

As is well known,  the Clifford algebra of the sum of two vector spaces is the ($\z_2$-graded) tensor product of the Clifford algebras of the two summands (This is denoted by  $\cl(V\oplus W)=\cl(V)\widehat\otimes\cl(W)$ in  \cite{LawMich}). Of course this is not restricted to two factors, but can be iterated. The case of two factors, i.e.\ $N=2$, will always be emphasized.

We consider the pullback bundles $p_i^*S_i\to M$ of $S_i$ over $M$ and once more we omit the projections in our further notation.  From $S_i$ we construct a bundle $W$ on $M$
\begin{align}
W&=\; (S_1^+\oplus S_1^-)\otimes\cdots \otimes (S_N^+\oplus S_N^-)\nonumber\\
&=\; \bigoplus_{\epsilon_1+\cdots+\epsilon_N=0} \bigotimes_{i=1}^{N}  S_i^{\epsilon_i} \oplus 
  \bigoplus_{\epsilon_1+\cdots+\epsilon_N=1} \bigotimes_{i=1}^{N}  S_i^{\epsilon_i}.
\label{W}
\end{align}
For $N=2$ this is
\begin{align*}
W=& (S_1^+\oplus S_1^-)\otimes (S_2^+\oplus S_2^-)\\
=&\big( S_1^+\otimes S_2^+\oplus S_1^-\otimes S_2^-\big)\oplus\big(S_1^+\otimes S_2^-\oplus S_1^-\otimes S_2^+\big)
\end{align*}
The bundles $S_i^{\epsilon_i}$  are subbundles of $S_i$ with the $\z_2$-degree defined by the label\footnote{We identify $\z_2=\{0=+,1=-\}$.}  $\epsilon_i=\z_2$. This induces on $W$ the  natural ${\z}_2$-grading given by the decomposition (\ref{W}). An element  $\Xi\in W$ is called totally homogenous, if $\Xi^\epsilon=\xi_1^{\epsilon_1}\otimes\cdots\otimes \xi_N^{\epsilon_n} $, i.e.\ it  has only contributions from one of the summands in the decomposition (\ref{W}). We will  use the multi index notation  $\epsilon=(\epsilon_1,\ldots,\epsilon_N)$.

The subbundles of $S_i$ are chosen in such a way that  the Clifford multiplication acts via $TM_i\otimes S_i^\pm\to S_i^\mp$, i.e.
\begin{equation}
X_i\xi_i=X_i{\binom{\xi^+_i}{\xi_i^-}}=\binom{X_i\xi^-_i}{X_i\xi^+_i}
   =\begin{pmatrix}0&X_i\\X_i&0\end{pmatrix}\xi_i.
\end{equation}
Clifford multiplication is an odd operation on $S_i^+\oplus S_i^-$, i.e.\ $|X_i\phi_i|=|\phi_i|+1$.

We will consider the continuation of this to the bundle $W$. Therefore we introduce the linear map $\delta_k :W\to W$, which is defined by  its action on totally homogenous elements $\Xi^\epsilon=\xi_1^{\epsilon_1}\otimes\cdots\otimes\xi_N^{\epsilon_N}$
\begin{equation}
\delta_k (\Xi^\epsilon)= (-)^{\epsilon_1+\cdots+\epsilon_{k-1}} \Xi^\epsilon  
\end{equation}
$\delta_k$ can be seen as the sign, we get by ``putting an odd operator acting on $S_k$ at the right place in the tensor product''. We define $\delta_1=\id$.

For $X=X_1+\cdots+ X_N\in TM$ and totally homogenous $\Xi^\epsilon=\xi_1^{\epsilon_1}\otimes\cdots\ \otimes \xi^{\epsilon_N}_N\in W$ we define
\begin{multline}
X\Xi^\epsilon\ :=\ \sum\nolimits_{i=1}^N ({\bf 1}\otimes\cdots \otimes  A_iX_i\otimes\cdots\otimes{\bf 1})\, \Xi^\epsilon\\
  = \sum\nolimits_{i=1}^N \delta_i(\Xi^\epsilon)\xi_1^{\epsilon_1}\otimes\cdots \otimes A_iX_i\xi_i^{\epsilon_i}\otimes\cdots\otimes \xi_N^{\epsilon_N}.
\end{multline}
In particular for $N=2$
\begin{multline*}
(X_1+X_2)\Xi^{(\epsilon_1,\epsilon_2)} = ({\bf 1} \otimes  A_1X_1+A_2X_2\otimes {\bf 1})\Xi^{(\epsilon_1,\epsilon_2)}\\ 
  =A_1X_1\xi_1^{\epsilon_1}\otimes\xi_2^{\epsilon_2}+(-)^{\epsilon_1}\xi_1^{\epsilon_1}\otimes A_2X_2\xi_2^{\epsilon_2}
 \end{multline*}

\begin{proposition}\label{prop1}
Let $(M,g)$, $(M_i,g_i)$, $S_i$ and $W$ be as before. 
The Clifford relation 
\begin{equation}
(XY+YX)\Xi =-2 g(X,Y)\Xi
\end{equation}
holds for all $X,Y\in TM$ and  $\Xi\in W$.
\end{proposition}

\begin{proof}
Because of linearity it is sufficient to prove the statement for $X_i\in TM_i$ and $Y_j\in TM_j$ and totally homogenous $\Xi^\epsilon$. We have to distinguish the cases $i=j$ and $i<j$. We recall the property  $\delta_k(X_i\Xi^\epsilon) =\begin{cases}-\delta_k(\Xi^\epsilon)&i<k\\\delta_k(\Xi^\epsilon)&i\geq k\end{cases}$. 

With these informations we get
\begin{align*}
(X_i& Y_j+Y_jX_i) \Xi^\epsilon = \\
=&\; X_i  \delta_j(\Xi^\epsilon)
    \xi_1^{\epsilon_1}\otimes\cdots \otimes A_jY_j\xi_j^{\epsilon_j}\otimes\cdots\otimes \xi_N^{\epsilon_N}\\
&  
+Y_j \delta_i(\Xi^\epsilon)
   \xi_1^{\epsilon_1}\otimes\cdots \otimes A_iX_i\xi_i^{\epsilon_i}\otimes\cdots\otimes \xi_N^{\epsilon_N}
\\
=&\; 
\begin{cases}
 \delta_i(A_jY_j\Xi^\epsilon)\delta_j(\Xi^\epsilon)
\xi_1^{\epsilon_1}\otimes\cdots\otimes A_iX_i\xi^{\epsilon_i}_i\otimes\cdots \otimes A_jY_j\xi_j^{\epsilon_j}
    \otimes\cdots\otimes \xi_N^{\epsilon_N} & \\[0.5ex]
 \  +\delta_j(A_iX_i\Xi^\epsilon) \delta_i(\Xi^\epsilon)
   \xi_1^{\epsilon_1}\otimes\cdots \otimes A_iX_i\xi_i^{\epsilon_i}\otimes\cdots\otimes A_jY_j\xi^{\epsilon_j}_j
   \otimes \cdots\otimes \xi_N^{\epsilon_N}
& i<j \\
 \delta_i(\Xi^\epsilon)\delta_i(\Xi^\epsilon)
 \xi_1^{\epsilon_1}\otimes\cdots \otimes (A_iX_i)(A_iY_i)\xi_i^{\epsilon_i}\otimes\cdots\otimes \xi_N^{\epsilon_N}
& \\
 \   +\delta_i(\Xi^\epsilon)\delta_i(\Xi^\epsilon)
   \xi_1^{\epsilon_1}\otimes\cdots \otimes (A_iY_i)(A_i X_i)\xi_i^{\epsilon_i}\otimes\cdots\otimes \xi_N^{\epsilon_N}
& i=j 
\end{cases}\\
=&\;
\begin{cases}
0 
&i<j\\[1ex]
\xi_1^{\epsilon_1}\otimes\cdots \otimes (A_i X_iA_iY_i+A_iY_iA_iX_i)\xi_i^{\epsilon_i}\otimes\cdots\otimes \xi_N^{\epsilon_N}=g(X_i,Y_i)\Xi^\epsilon
&i=j
\end{cases}\\[1ex]
=&\; -2g_i(A_iX_i,A_iY_j)\Xi^\epsilon\; =\; -2g(X_i,Y_i)\Xi^\epsilon
\end{align*}
\end{proof}
Up to now we did not specify the subbundles of $S_i^{\epsilon_i}\subset S_i$. We have to distinguish two different situations:

\begin{itemize}
\item If $\dim M_i=2n_i$ is even, the bundle $S_i$ (of rank $2^{n_i}$)  itself admits a natural $\z_2$-grading induced by the volume element. And we take exactly this one.
\item If $\dim M_i=2n_i+1$ is odd, the spinor bundle $S_i$ (of rank $2^{n_i}$) does not admit a natural $\z_2$-grading. In this case we double the bundle and define $S_i^+:=S_i$ and $S_i^-:=\Pi_iS_i$.
\end{itemize}

In this definition of the subbundles, $\Pi_i: S_i^+\oplus S_i^-\to S_i^-\oplus S_i^+$ denotes the parity operator 
with $(\Pi_i(S_i^+\oplus S_i^-))^\pm=S_i^\mp$  and $\Pi_i^2=\id$. Explicitly we have
\begin{equation}
\Pi_i\binom{\xi_i^+}{\xi_i^-}=\binom{\xi_i^-}{\xi_i^+}.
\end{equation}
The parity operator is naturally extended to $W$ by its action on totally  homogenous elements 
\begin{equation}
\Pi (\Xi) 
=\frac{1}{\sqrt{N}}\sum_{i=1}^{N} (-)^{\epsilon_1+\cdots+\epsilon_{i-1}}
\xi_1^{\epsilon_1} \otimes\cdots\otimes \Pi_i\xi_i^{\epsilon_i}\otimes \cdots\otimes \xi_N^{\epsilon_N}
\end{equation}
where we use that the parity operator is formally odd. The proof for $\Pi^2=\id$ is similar to that for the Clifford relation.

Let $N_e$ and $N_o$ be the number of even and odd  dimensional manifolds in the product $M=M_1\times\cdots\times M_N$, respectively (i.e.\ $N_e+N_o=N$). The dimension of the product manifold is $\dim M=2n+N_o$  with $n:=\sum_{i=1}^N n_i$ and the rank of the bundle $W$ is $2^{n+N_o}$. 

The spinor bundle of $M$ should have the rank $2^{n+[\frac{N_o}{2}]}$. We will construct a subbundle $S$  of $W$ of this rank.

The trick is to diagonalize some of the bundles  which come from the odd dimensional manifolds. The diagonalization is denoted by $\Delta$ and the bundle $S$  is constructed as follows:

Choose $N_o-[\frac{N_o}{2}]$ of the odd dimensional bundles\footnote{This does not depend on the choice, because the resulting (non graded) bundles are isomorphic.} and consider the subbundle $S\subset W$ given by 
\begin{equation}
S:= \bigotimes_{i=1}^{N_e} S_i 
\otimes \bigotimes_{j=N_e+1}^{N- N_o+[\frac{N_o}{2}]}( S_j\oplus \Pi_j S_j)
\otimes\bigotimes_{k=N-N_o+[\frac{N_o}{2}]+1}^{N}\Delta(S_k\oplus \Pi_k S_k)
\label{spinorbundleproduct}
\end{equation}
which has the rank $2^K$ with
\begin{align}
K&= {\sum_{i=1}^{N_e} n_i   +\sum_{j=N_e+1}^{N-N_o+[\frac{N_o}{2}]}n_j+
    (N-N_o+[\frac{N_o}{2}]-N_e)  +\sum_{k=N-N_o+[\frac{N_o}{2}]+1}^{N}n_k}\nonumber\\
&= \sum_{i=1}^{N} n_i+[\frac{N_o}{2}]\ =\ n+[\frac{N_o}{2}]\nonumber
\end{align}

We specialize this to the case $N=2$:

\begin{itemize}
\item[$\bullet_1$] $\dim M_1=2n$, $\dim M_2=2m$:$\quad$  
$S$ over  $M=M_1\times M_2$ is given by 
\begin{equation}
S=S_1\otimes S_2
\end{equation}
\item[$\bullet_2$] $\dim M_1=2n$, $\dim M_2=2m+1$:$\quad$  
The subbundle $S$ over $M$ is given by  
\begin{equation}
S = S_1\otimes \Delta(S_2\oplus \Pi_2 S_2)\simeq S_1\otimes S_2.
\end{equation}
\item[$\bullet_3$] $\dim M_1=2n+1$, $\dim M_2=2m+1$:$\quad$  
In this case $S$ is defined by
\begin{multline*}
S=(S_1\oplus \Pi_1 S_1)\otimes\Delta(S_2\oplus \Pi_2 S_2)\simeq (S_1\oplus S_1)\otimes S_2\\
 \simeq S_1\otimes S_2\oplus S_1\otimes S_2\simeq (S_2\oplus S_2)\otimes S_1,
\end{multline*}
where we have to emphasize that the equivalences may not respect the $\z_2$-grading.
\end{itemize}

That this bundle is indeed a spinor  bundle will be clear from the following remark.
Although we have established the Clifford multiplication, the bundle is not a priori a spinor bundle.
For that it should be constructed as an associated vector bundle to the $ Spin$-principal bundle $P_{ Spin(D)}(M)$ or of one reduction of this. If such a reduction does not exist in an appropriate way, then $W$ is not of this kind. The reason is that we do not have an action of $ Spin(D)$ on the standard fibre of $W$.

The maximal subgroup of $ Spin(D)$ which is compatible with the structure of the standard fibre is
\begin{multline}
S( Pin(D_1)\times\cdots\times Pin(D_N)):= Pin(D_1)\times\cdots\times Pin(D_N)\cap\cl^+(D)\\
\subset\otimes_{i=1}^N\cl(D_i).
\end{multline}
This contains the subgroup $ Spin(D_1)\times\cdots\times Spin(D_N)$, which will be of interest soon.

In the next proposition we show that we are indeed able to write $W$ as an associated bundle, if we demand the reduction of the structure group of $M$ to $SO(D_1)\times\cdots\times SO(D_N)$ compatible with the natural splitting (\ref{tangentsplitting}). This is a weaker condition than (\ref{metricsplitting}), which was at least necessary to get  the Clifford multiplication. 

\begin{proposition}\label{prop2}
Let $M$ be a Riemannian spin manifold. If the structure group of $M$ can be reduced to  $SO(D_1)\times\cdots\times SO(D_N)$, we have a reduction of the spin principal bundle to $ Spin(D_1)\times\cdots\times Spin(D_N)$.
\end{proposition}

\begin{proof}
We use the following notations: $G:= Spin(D)$, $\widetilde G:= Spin(D_1)\times\cdots\times Spin(D_N)$, $H:=SO(D)$ and $\widetilde H:= SO(D_1)\times\cdots\times SO(D_N)$.

Let $P_H$ be the principal bundle of orthonormal frames of $(M,g)$ and  $P_{\widetilde H}$ be the reduction to $\widetilde H$, which we denote by $\imath$.  Furthermore let $P_G$ the $Spin$-principal bundle over $M$ with the two fold covering $\lambda :P_G\to P_H$, compatible with the right action of $G$ and $H$ and the cover $\lambda:G\to H$ for which we take the same symbol.  

We collect this by writing $P_G\stackrel{\lambda}{\longrightarrow}P_H\stackrel{\imath}{\longleftarrow} P_{\widetilde H}$. 
We use the pullback construction for fibre bundles cf.\ (\cite{MichorBook}) to complete this edge to a commutative diagram and we denote the pullback by  $(P_G\times P_{\widetilde H}) \diagup P_H$. 
\begin{equation}
{
\xymatrix{
P_G  \ar[r]^-{\lambda}                                                        & P_H   \\
(P_G\times P_{\widetilde H}) \diagup {P_H}\ar[u]\ar[r]  & P_{\widetilde H}\ar[u]_-{\imath}
}}
\end{equation}
Its total space is given by
\begin{equation}
\left\{
(p_G,p_{\widetilde H})\in P_G\times P_{\widetilde H}\,|\,\lambda(p_G)=\imath(p_{\widetilde H})
\right\}
\end{equation}
The bundle $(P_G\times P_{\widetilde H}) \diagup {P_H}$ is a principal bundle with the right action of its standard fibre  
$G\times\widetilde H\diagup H=\{(g,\tilde h)\,|\,\lambda(g)=\imath(\tilde h)=\tilde h\}
\simeq \lambda^{-1}(\widetilde H)\simeq \widetilde G$ 
on the total space  defined in the obvious way. We take the action on the Cartesian product 
$(p_G,p_{\widetilde H}) a=(p_G a ,p_{\widetilde H}\lambda(a))$ which is compatible with the quotient:
\begin{equation}
\lambda( p_G a)=\lambda(p_G)\lambda a =\imath(p_{\widetilde H})\lambda a =\imath(p_{\widetilde H}\lambda(a)) 
\end{equation}
This shows that the $Spin$-principal bundle $P_G$ is reduced to a $Spin(D_1)\times\cdots\times Spin(D_N)$-principal bundle. The reduction is just the left vertical arrow in the diagram.
\end{proof}

From this we get
\begin{theorem}\label{prop3}
Let  $(M=M_1\times\cdots\times M_N,g)$ be a Riemannian  manifold, given as a (not necessarily Riemannian) product of  the simply connected Riemannian spin manifolds  $(M_i,g_i)$. The metric $g$ is connected to the metrics $g_i$ via (\ref{metricsplitting}).

Then $M$ is spin and the spinor bundle is isomorphic to the subbundle $S$ of $W$ constructed in this section. 
\end{theorem}
\begin{proof}
That $M$ is spin, follows immediately from the fact that the second Stiefel-Whitney class behaves additively under products of manifolds.  The ON-frame bundles w.r.t.\ $g$ and $g_1+\cdots+g_N$  are isomorphic cf.\ (\ref{isomorphism}). 
On $M$ there exists exactly one spin structure, because $M$ is simply connected. 
So the spin principal bundles obtained by $g$ and $g_1+\cdots +g_N$ are isomorphic. 
From the last proposition we get that the $Spin$-principal bundle $P_G$ over $(M,g)$  is reducible to $Spin(D_1)\times\cdots\times Spin(D_N)$. 

With  the notations from the previous proposition and the construction for $S\subset W$ we have established the following isomorphism 
\begin{equation}
S\simeq \big((P_G\times P_{\widetilde H}) \diagup {P_H}\big)\times_{\widetilde G} \hat S
\end{equation}
where $\hat S$ denotes the standard fibre of $S$.
\end{proof}

\begin{remark}
An important example of this construction is the case $N=2$ with one of the factors being one dimensional and the metric is a warped product. This has been discussed in detail in  \cite{Baum5} and \cite{Baum4}. \label{Baum}
\end{remark}

From the construction it is clear that the spin connection obtained from the Levi-Civita connection of $g$ need not be compatible with the tensor structure of $S$. Explicit formulas for the connection in the case of an one dimensional factor and warped products can also be found in \cite{Baum4}.  By claiming that the holonomy of $M$ is contained in $SO(D_1)\times\cdots\times SO(D_N)$ -- which is the same as to say that the projections $p_i$ are parallel -- we make sure this further compatibility.
As we mentioned above, this further assumption forces $M$ to be a local Riemannian product (in the case of $M$ being simply connected and complete this decomposition is global).

This constructions yields for $A=id$ the following
\begin{proposition}\label{prop4}
Let $M$ be the Riemannian product of $(M_i,g_i)$. Then the spinor bundle of $M$ with respect to the induced spin structure is given by $S$ from the construction above.
\end{proposition}

We add some remarks, which explain our further notation and draw the attention to the pseudo-Riemannian case.
\begin{remark}
\begin{enumerate}
\item Proposition \ref{prop4} is also true in the case of metrics which are not of Euclidean signature. 
\item For the special case $N=D$, i.e.\ $D_i=1$, our constructions ends up with  the Clifford action on  $\c^2\otimes\cdots\otimes \c^2$ as given in \cite{BFGK}.
\item For another special case $N=2$ and $D_2=1$ this construction leads to the discussion in \cite{Baum5}.
\item For $N=2$ in the cases $\bullet_1$ and $\bullet_2$, our construction yields the -- at least in physics literature --  frequently used decomposition of the $\gamma$-matrices $\Gamma^A$, for $1\leq A\leq D$ in tensor products of $\gamma$-matrices of the respective factors $\gamma^a$, $\gamma^\alpha$ for $1\leq a\leq D-2n=D_2$, $1\leq\alpha\leq 2n=D_1$. The $\z_2$-grading is ensured by using the volume element $\hat \gamma:=\gamma^1\cdots\gamma^{2n}$ which anticommutes with all $\gamma^a$  and define 
\begin{equation}
\Gamma^a=\gamma^a\otimes\hat\gamma, \qquad \Gamma^\alpha={\bf 1}\otimes\gamma^\alpha\;,
\label{gammasplitting}
\end{equation}
compare e.g.\ \cite{DuffNilPop1}.
\end{enumerate}
\end{remark}


\end{document}